\documentclass[10pt]{amsart}

\usepackage[latin1]{inputenc}
\usepackage{mathptmx}
\usepackage{amsmath}
\usepackage{dsfont}
\usepackage{amsfonts}
\usepackage{amssymb}
\usepackage{amsthm}
\usepackage{mathrsfs}
\usepackage{MnSymbol}
\usepackage{color,soul}
\usepackage[text={4.5in,7.2in},centering]{geometry}
\usepackage{gensymb}
\usepackage{caption}
\usepackage{float}
\theoremstyle{definition}
 
\numberwithin{dummy}{section}

\theoremstyle{plain}
\newtheorem{thm}{Theorem}[section]

\newtheorem{prop}[thm]{Proposition}

\theoremstyle{definition}
\newtheorem{defn}{Definition}[section]

\newtheorem{exmp}{Example}[section]

\theoremstyle{remark}

\DeclareMathAlphabet\mathbb{U}{msb}{m}{n}

\let\oldhat\hat                

\renewcommand{\hat}[1]{\oldhat{\mathbf{#1}}}

\usepackage{calc}
\usepackage{tikz}

\usetikzlibrary{decorations.markings}
\usetikzlibrary{positioning,automata}
\usetikzlibrary{shapes}
\usetikzlibrary{arrows}
\usetikzlibrary{fit}
\usepackage{verbatim}
\tikzset{
	>=stealth',
	vertex/.style={
		circle,
		draw=#1,
		fill=#1,
		inner sep = 2pt,
		outer sep = 0pt,
		text centered},
	edge/.style={
		-,
		thin,
		draw=black,
		},
	text style/.style={
		sloped,
		text=black,
		font=\normalsize,
		above}
	}

\tikzset{every loop/.style={min distance=10mm,in=45,out=135,looseness=0}}

\title{$P$-Matchings in Graphs: A Brief Survey with Some Open Problems}

\author{Todd Fenstermacher, Soumendra Ganguly, Stephen Hedetniemi, \\ Renu Laskar \\ Clemson University}

\begin{document}
	
	\maketitle
		
	\begin{abstract}
		 For a graph $G=(V,E),$ a matching $M$ is a set of independent edges. The topic of matchings is well studied in graph theory. In this paper many varieties of matchings are discussed.
	\end{abstract}

\section{Introduction}

We begin with some standard definitions on a graph $G = (V,E).$
For a vertex $v \in V(G),$ the \textit{open neighborhood} of $v$ is $N(v) = \{ u | uv \in E(G)\},$ and the \textit{closed neighborhood} of $v$ is $N[v] = N(v) \cup \{v\}.$ A set $S \subseteq V(G)$ is a \textit{dominating set} if for all $u \in V(G) - S,$ there exists $v \in S$ such that $uv \in E(G).$ The domination number of $G$ is $\gamma(G) = \min \{ |S| : S \text{ is dominating set}\}.$ A set $S \subseteq V(G)$ is \textit{independent} if $u,v \in S$ implies $uv \notin E(G).$ The independence number is $\beta_0(G) = \max \{ |S| : S \text{ is independent} \}.$ A set $S \subseteq V(G)$ is a \textit{vertex cover} if every edge of $G$ is incident to a vertex in $S.$ The vertex covering number of $G$ is $\alpha_0(G) = \min\{ |S| : S \text{ is a vertex cover}\}.$ A set $F \subseteq E(G)$ is an \textit{edge cover} if every vertex of $G$ is incident with an edge in $F.$ The edge covering number of $G$ is $\alpha_1(G) = \min\{ |F| : F \text{ is an edge cover}\}.$ A set $M \subseteq E(G)$ is a \textit{matching} if no two edges in $M$ have a vertex in common. The matching number of $G$ is $\beta_1(G) = \max\{ |M| : M \text{ is a matching}\}.$ If a vertex $u$ is incident with some edge in $M,$ then $u$ is said to be saturated by $M.$ A matching $M$ is \textit{maximal} if for all $e \notin M,$ $M\cup\{e\}$ is not a matching. The lower matching number of $G$ is $\beta^-_1(G) = \min\{ |M| : M \text{ is a maximal matching}\}.$ The matching number and lower matching number are demonstrated in Figure \ref{P8}. A matching $M$ is  \textit{perfect} if it saturates all the vertices of $G.$ Not every graph has a perfect matching.  

\begin{center}
	\begin{figure}[H]
	\begin{tikzpicture}
	\node[vertex, label = 1] (0) {};
	\node[vertex, right = .5cm of 0, label = 2] (1) {}
	edge[edge] node {} (0);	
	\node[vertex, right= .5cm of 1, label = 3] (2) {}
	edge[edge] node {} (1);
	\node[vertex, right= .5cm of 2, label = 4] (3) {}
	edge[edge] node {} (2);
	\node[vertex, right= .5cm of 3, label = 5] (4) {}
	edge[edge] node {} (3);
	\node[vertex, right= .5cm of 4, label = 6] (5) {}
	edge[edge] node {} (4);
	\node[vertex, right = .5cm of 5, label = 7] (6) {}
	edge[edge] node {} (5);	
	\node[vertex, right= .5cm of 6, label = 8] (7) {}
	edge[edge] node {} (6);			
	\end{tikzpicture}
	\caption{$\beta_1(P_8) = 4  \text{ from }	\{ 12, 34, 56, 78	\}$ and $\beta^{-}_1(P_8) = 3 \text{ from } \{23, 45, 67  \}$}
	\label{P8}
\end{figure}
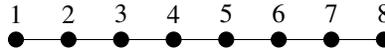
\end{center}

\section{Historical Background}

Here we present some of what are perhaps the most significant results of matching theory. These results can be found in any standard textbook on graph theory.  An excellent book on matching theory is \textit{Matching Theory} by Lov\'{a}sz and Plummer \cite{Lovasz}. In particular, the preface of this book gives an extraordinary history of some of these results. 

\begin{thm}[Frobenius (1917): Marriage Theorem ]\label{Marriage}
	$G = (A,B,E)$ bipartite has a perfect matching if and only if
	
	i) $|A| = |B|$
	
	ii) $|X| \leq |N(X)|$ for each $X \subseteq A$
\end{thm}

\begin{thm}[K\"{o}nig-Egerv\'{a}ry (1931)]\label{Konig}
	$G$ bipartite implies $\beta_1(G) = \alpha_0(G).$
\end{thm}

\begin{defn}
	Let $A = \{A_1,A_2,\dots, A_n\}$ be a collection of subsets of a set $X \subseteq V(G).$ A \textit{system of distinct representatives} (SDR) for $A$ is a set of distinct elements $a_1,a_2,\dots,a_n$ in $X$ such that $a_i \in A_i.$
\end{defn}

\begin{thm}[P. Hall (1935) : ``Hall's Condition"] \label{Hall}
	$A = \{A_1, A_2, \dots A_m\}$ a set of subsets of a set $X$ has a SDR $\iff$ $| \cup_{i \in S} A_i| \geq |S|$ for every $S \subseteq \{1,2,\dots,m\}.$ 
\end{thm}

Note that if $|X| = |A|,$ then Theorem \ref{Hall} is equivalent to Theorem \ref{Marriage}. A great book which explores problems having their origin in above theorem by Hall is \textit{Transversal theory} by Mirsky \cite{Mirsky}.
We also note that a number of early equivalent results were proved by several authors. Historically the first of these results was by Frobenius, but the names of K\"{o}nig, Egerv\'{a}ry, and Hall are more often assigned to various generalizations of this result. In particular,  Theorem \ref{Konig} is an important example of a class of results known as ``minimax theorems." For example, the well-known max-flow min-cut theorem is one such minimax theorem.

The last of our classic theorems in this section relates four of the parameters discussed in the introduction. 

\begin{thm}[Gallai Identities (1959)]
	If $G$ has no isolates, then
	
	i) $\alpha_0(G) + \beta_0(G) = |V(G)|$
	
	ii) $\alpha_1(G) + \beta_1(G) = |V(G)|$
\end{thm}

Theorems on maximum and perfect matchings came later, and attempted to answer such questions as, ``When is a matching maximum?" or, ``When does a perfect matching exist?" Berge and Tutte  and many others have attempted to answer such questions \cite{Berge1},\cite{Berge2},\cite{Tutte}.

Lastly we note that matching theory has not developed in a vacuum, but has had an influence in other areas of mathematics. For example the minimax theorems as mentioned earlier are close to the birth of duality theory in linear programming. Furthermore, the first nontrivial polyhedron studied by Edmonds was the matching polytope, and this lead to ground-breaking work in the areas of facet determination and ``good" characterizations. 
Even more so, Edmond's matching algorithm had influence in complexity theory by showing that the challenging problem of matching was solvable in polynomial time. 

\section{Varieties of Matchings}

In this section we highlight a few interesting varieties of matchings. The varieties discussed here are examples of what we call $P$-matchings.

\begin{defn}
Let $P$ be a property of $G,$ and let $M$ a matching of $G$ such that the induced subgraph of $M,$ $\langle M \rangle,$ has property $P,$ then we call $M$ a $P$-matching. 
\end{defn}

With regard to $P$-matchings, we are interested in two parameters: 
\[ \beta_P = \max \{ |M| : \langle M \rangle \text{ has property } P  \}\]
\[\beta^{-}_P = \min \{ |M| : \langle M \rangle \text{ is maximal with respect to } P  \}.\]

\begin{exmp}[Induced Matching (strong matching)] 
$M$ is an induced (strong) matching if $\langle M \rangle$ is a disjoint union of $K_2$'s. 

%

Strong matching was introduced by Cameron in 1989 \cite{Cameron}. The two parameters of interest here are $\beta_*(G) = \max \{ |M| : M \text{ is induced matching} \}$ and $\beta^-_*(G) = \min \{ |M|: M \text{ is maximal induced matching}\}.$
\end{exmp}

\begin{figure}[H]
	\begin{tikzpicture}
	\node[vertex ] (0) {};
	\node[vertex, right = 1cm of 0] (1) {};
	\node[vertex, right= 1cm of 1] (2) {};
	\node[vertex, right= 1cm of 2] (3) {};
	\node[vertex, below= .5cm of 0] (4) {};
	\node[vertex, below= .5cm of 1] (5) {}
	edge[edge] node {} (0);
	\node[vertex, below = .5cm of 2] (6) {}
	edge[edge] node {} (0)
	edge[edge] node {} (1)
	edge[edge] node {} (3);	
	\node[vertex, below= .5cm of 3] (7) {};			
	
	\node[vertex, right = 2cm of 3] (8) {};
	\node[vertex, right = 1cm of 8] (9) {};
	\node[vertex, right= 1cm of 9] (10) {};
	\node[vertex, right= 1cm of 10] (11) {};
	\node[vertex, below= .5cm of 8] (12) {}
	edge[edge] node {} (9);
	\node[vertex, below= .5cm of 9] (13) {}
	edge[edge] node {} (10);
	\node[vertex, below = .5cm of 10] (14) {}
	edge[edge, red] node {} (10)
	edge[edge] node {} (11);
	\node[vertex, below= .5cm of 11] (15) {}
	edge[edge] node {} (8);
	
	\path[draw,-,ultra thick] (0) -- (4);
	\path[draw,-,ultra thick] (1) -- (5);
	\path[draw,-,ultra thick] (2) -- (6);
	\path[draw,-,ultra thick] (3) -- (7);
	
	\path[draw,-,ultra thick] (8) -- (12);
	\path[draw,-,ultra thick] (9) -- (13);
	\path[draw,-,ultra thick] (10) -- (14);
	\path[draw,-,ultra thick] (11) -- (15);

	\end{tikzpicture} 
	\caption{Examples of ur (left) and not ur (right) matchings.}
	\label{ur}
	\end{figure}
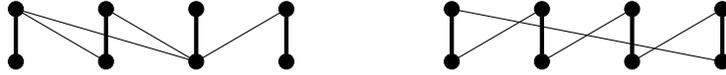

\begin{exmp}[Uniquely Restricted Matching (ur)]
A matching $M$ is a uniquely restricted if the only perfect matching of $\langle M \rangle$ is $M$. Uniquely restricted matching was introduced by Golumbic et al. in 2001 \cite{Golumbic1}. The parameters of interest are\[ \beta_{ur}(G) = \max \{ |M| : M \text{ is a ur matching}\}\] and \[ \beta^-_{ur}(G) = \min \{ |M| : M \text{ is a maximal ur matching}\}.\]
\end{exmp}
Figure \ref{ur} contains an example of a matching which is ur and a matching which is not ur. 

The next four examples of of $P$-matchings are \textit{connected} matchings, \textit{isolate free} matchings, \textit{disconnected matchings}, and \textit{acyclic} matchings. These varieties were introduced by Goddard et al. in 2005 \cite{Goddard}. 

\begin{exmp}[Connected matching]
	A matching $M$ is connected if $\langle M \rangle$ is connected. We have $\beta_c(G) = \max\{ |M| : M \text{ is connected matching} \}$ and $\beta^-_c(G) = \min\{ |M| : M \text{ is maximal connected matching}\}.$ 
\end{exmp}

\begin{exmp}[Isolate free matching]
	A matching $M$ is isolate free if $|M| = 1$ or $\langle M \rangle$ has no $K_2$ component. We have $\beta_{if}(G) = \max\{ |M| : M \text{ is isolate free matching} \}$ and $\beta^-_{if}(G) = \min\{ |M| : M \text{ is maximal isolate free matching}\}.$ 
\end{exmp}

\begin{exmp}[Disconnected matching]
	A matching $M$ is disconnected if $|M| = 1$ or $\langle M \rangle$ is disconnected. We have $\beta_{dc}(G) = \max\{ |M| : M \text{ is disconnected matching} \}$ and $\beta^-_{dc}(G) = \min\{ |M| : M \text{ is maximal disconnected matching}\}.$ 
\end{exmp}

\begin{exmp}[Acyclic matching]
	A matching $M$ is acyclic if $\langle M \rangle$ is acyclic. We have $$\beta_{ac}(G) = \max\{ |M| : M \text{ is acyclic matching} \}$$ and $$\beta^-_{ac}(G) = \min\{ |M| : M \text{ is maximal acyclic matching}\}.$$ 
\end{exmp}

We next give a few simple results concerning these new varieties of matchings. 
\begin{prop}
	Let $G$ be a graph, then 
	
	i) $\beta_*(G) \leq \beta_{ac}(G) \leq \beta_{ur}(G) \leq \beta_1(G)$
	
	ii) $\beta_*(G) \leq \beta_{dc}(G) \leq \beta_1(G)$
	
	iii) $\beta_c(G) \leq \beta_{if}(G) \leq \beta_1(G)$	
\end{prop}

The following theorem is found in \cite{Golumbic1}. 
\begin{thm}
	A matching $M$ is ur if and only if $\langle M \rangle$ does not contain an alternating cycle with respect to $M.$
\end{thm}

This next theorem is found in \cite{Goddard}. 
\begin{thm}
	If $G$ is connected, then $\beta_c(G) = \beta_{if}(G) = \beta_1(G).$	
\end{thm}

The next variety of matching was introduced by Nordhaus in 1977 \cite{Alavi}. 
\begin{exmp}
	A set $M \subseteq V \cup E$ is a \textit{total matching} of $G$ if the elements of $M$ are pairwise independent, and $M$ is maximal.
\end{exmp}

%

In a 1970 paper \cite{Graham}, Graham introduced what he called a simple cutset. This gave the basis for what we call a separating matching. 
\begin{exmp}
	A matching $M$ is \textit{separating} (or \textit{disconnecting}) matching if $M$ is an edge cut.
\end{exmp}
 
The hypercube $Q_n$ is a classic example of a graph with a separating matching. 

%
%
%
%
%
%

The last variety of matchings that we mention here before introducing new varieties of matchings is $b$-matchings. We first need some notation. We let $d(v)$ denote the degree of a vertex $v,$ and let $b(v)$ be bound such that $0 \leq b(v) \leq d(v).$  
\begin{defn}
	A set $M \subseteq E(G)$ is a \textit{b-matching} of $G$ if $|\{(u,v) \in M\}| \leq b(v)$ for all $v.$ 	
\end{defn}

Much study has been done on the complexity of finding $b$-matchings, in particular Goodman, Hedetniemi, and Tarjan found a linear time greedy algorithm for trees \cite{Goodman}.

\section{New Classes of Matchings}

In this section we introduce what we believe to be new classes of matchings. We start with some definitions. Let $S \subseteq V(G)$ and let $u \in S.$ We say $u$ is \textit{irredundant} with respect to $S$ if $u$ has a private neighbor, i.e., $N[u] - N[S-\{u\}] \neq \emptyset.$ A vertex $v$ is an \textit{external private} neighbor of $u$ if $v$ is a private neighbor of $u$ and $v \notin S.$ A set $S$ is \textit{irredundant} if every vertex of $S$ is irredundant. We now introduce vertex-irredundant and edge-irredundant matchings. 

\begin{defn}
	A matching $M$ is vertex-irredundant if for every $e = uv \in M,$ either $u$ or $v$ has an external private neighbor.
\end{defn}

\begin{defn}
	A matching $M$ is edge-irredundant if every $e \in M$ has an edge $e' \notin M$ incident to $e$ and no other edges of $M.$
\end{defn}

For the above classes, we are interested in the usual parameters but introduce the following notation. 
\[ \beta^v_{IR}(G) = \max\{ |M| : M \text{ is vertex-irredundant matching}\}	\]
\[ \beta^v_{ir}(G) = \min\{ |M| : M \text{ is maximal vertex-irredundant matching}\}	\]
\[ \beta^e_{IR}(G) = \max\{ |M| : M \text{ is edge-irredundant matching}\}	\]
\[ \beta^e_{ir}(G) = \min\{ |M| : M \text{ is maximal edge-irredundant matching}\}	\]

We next introduce \textit{independent} and \textit{bipartite} matching. 

\begin{defn} 
	A matching is \textit{independent} if it has an orientation $(X,Y)$ such that $X$ is independent. 
\end{defn}

\begin{figure}[H]
	\begin{tikzpicture}
	
	\node[vertex, label = 1 ] (1) {};
	\node[vertex, right = 1cm of 1, label = 2] (2) {};
	\node[vertex, right= 1cm of 2, label = 3] (3) {}
	edge[edge] node {} (2);
	\node[vertex, right= 1cm of 3, label = 4] (4) {}
	edge[edge] node {} (3);
	
	\node[vertex, below= .5cm of 1, label = {below: 1'}] (1') {}
	edge[edge] node {} (1);
	\node[vertex, below= .5cm of 2,label = {below: 2'}] (2') {}
	edge[edge] node {} (1')
	edge[edge] node {} (2);
	\node[vertex, below = .5cm of 3,label = {below: 3'}] (3') {}
	edge[edge] node {} (3);	
	\node[vertex, below= .5cm of 4,label = {below: 4'}] (4') {}
	edge[edge] node {} (4);			
	
	\node[below right = .5cm and 1cm of 4, label=$\rightarrow$] (5) {};
	
	\node[right = 1.5cm of 4, label = {right: $X:$}] (0) {};
	\node[below = .5cm of 0, label = {right: $Y:$}] (0') {};
	\node[vertex, right = 1cm of 0, label = 1] (11) {};
	\node[vertex, right = 1cm of 11, label = 2] (12) {};
	\node[vertex, right= 1cm of 12, label = 3'] (13') {};
	\node[vertex, right= 1cm of 13', label = 4] (14) {};
	
	\node[vertex, below= .5cm of 11,label = {below: 1'}] (11') {}
	edge[edge] node {} (11);
	\node[vertex, below= .5cm of 12,label = {below: 2'}] (12') {}
	edge[edge] node {} (11')
	edge[edge] node {} (12);
	\node[vertex, below = .5cm of 13',label = {below: 3}] (13) {}
	edge[edge] node {} (12)
	edge[edge] node {} (13')
	edge[edge] node {} (14);
	\node[vertex, below= .5cm of 14,label = {below: 4'}] (14') {}
	edge[edge] node {} (14);			
	\end{tikzpicture}
	\caption{An independent matching}
	\label{Ind}	
\end{figure}
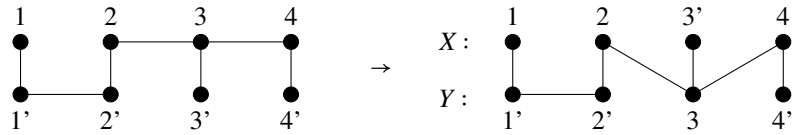

An example of an independent matching is given in Figure \ref{Ind}. 

%

The two parameters of interest are 
\[\beta_i(G) = \max \{ |M| : M \text{ is independent matching}\}\]
\[\beta^-_i(G) = \min \{ |M| : M \text{ is maximal independent matching}\}.\]

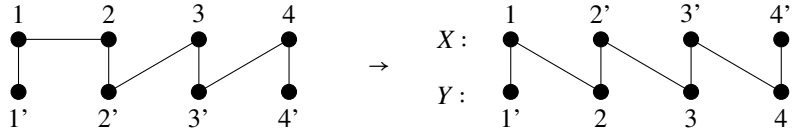
\begin{figure}[H]
	\begin{tikzpicture}
	\node[vertex, label = 1 ] (1) {};
	\node[vertex, right = 1cm of 1, label = 2] (2) {}
	edge[edge] node {} (1);
	\node[vertex, right= 1cm of 2, label = 3] (3) {};
	\node[vertex, right= 1cm of 3, label = 4] (4) {};
	
	\node[vertex, below= .5cm of 1, label = {below: 1'}] (1') {}
	edge[edge] node {} (1);
	\node[vertex, below= .5cm of 2,label = {below: 2'}] (2') {}
	edge[edge] node {} (3)
	edge[edge] node {} (2);
	\node[vertex, below = .5cm of 3,label = {below: 3'}] (3') {}
	edge[edge] node {} (4)
	edge[edge] node {} (3);	
	\node[vertex, below= .5cm of 4,label = {below: 4'}] (4') {}
	edge[edge] node {} (4);			
	
	\node[below right = .5cm and 1cm of 4, label=$\rightarrow$] (5) {};
	
	\node[right = 1.5cm of 4, label = {right: $X:$}] (0) {};
	\node[below = .5cm of 0, label = {right: $Y:$}] (0') {};			
	\node[vertex, right = 1cm of 0, label = 1] (11) {};
	\node[vertex, right = 1cm of 11, label = 2'] (12') {};
	\node[vertex, right= 1cm of 12, label = 3'] (13') {};
	\node[vertex, right= 1cm of 13', label = 4'] (14') {};
	
	\node[vertex, below= .5cm of 11,label = {below: 1'}] (11') {}
	edge[edge] node {} (11);
	\node[vertex, below= .5cm of 12,label = {below: 2}] (12) {}
	edge[edge] node {} (11)
	edge[edge] node {} (12');
	\node[vertex, below = .5cm of 13',label = {below: 3}] (13) {}
	edge[edge] node {} (12')
	edge[edge] node {} (13');
	\node[vertex, below= .5cm of 14,label = {below: 4}] (14) {}
	edge[edge] node {} (13')
	edge[edge] node {} (14');			
	\end{tikzpicture}
	\caption{A bipartite matching}
	\label{Bi}	
\end{figure}

\begin{defn}
	A matching is \textit{bipartite} if it has an orientation $(X,Y)$ such that both $X$ and $Y$ are independent. 
\end{defn}
An example of a bipartite matching is given in Figure \ref{Bi}. 

The parameters of interest for a bipartite matching are 
\[\beta_b(G) = \max \{ |M| : M \text{ is bipartite matching}\}\]
\[ \beta^-_b(G) = \min \{ |M| : M \text{ is maximal bipartite matching}\}.\]

We next introduce two more related varieties of matchings. But first we need the following definitions. 
\begin{defn}
	Two edges $e_1, e_2 \in E(G)$ are \textit{closed neighborhood adjacent (cnbr)} if there exists $v \in V(G)$ such that $e_1e_2 \in E(\langle N[v]\rangle).$ We say that $e_1, e_2$ are \textit{cnbr independent} if they are not cnbr adjacent.  
\end{defn}

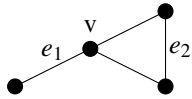
\begin{figure}[H]
	\begin{tikzpicture}
	\draw node at (0,0) [vertex, label = v] [] {};
	\draw node at (-1,-.5) [vertex] [] {};
	\draw node at (1,.5) [vertex] [] {};
	\draw node at (1,-.5) [vertex] [] {};
	
	\draw[-,decorate] (0,0) -- (-1,-.5) node [midway, yshift=.2cm ] {$e_1$};
	\draw[-] (0,0) -- (1,.5);
	\draw[-] (0,0) -- (1,-.5);
	\draw[-,decorate] (1,.5) -- (1,-.5) node [midway, xshift =.2cm] {$e_2$};

	\end{tikzpicture}
	\caption{$e_1$ and $e_2$ are cnbr adjacent}
	\label{cnbr}
\end{figure}

Figure \ref{cnbr} gives an example of two edges which are cnbr adjacent. 

\begin{defn}
	Two edges $e_1, e_2 \in E(G)$ are \textit{open neighborhood adjacent (onbr)} if there exists $v \in V(G)$ such that $e_1e_2 \in E(\langle N(v)\rangle).$ We say that $e_1, e_2$ are \textit{onbr independent} if they are not onbr adjacent. 
\end{defn}

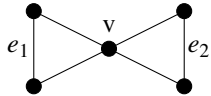
\begin{figure}[H]
	\begin{tikzpicture}
	\draw node at (0,0) [vertex, label = v] [] {};
	\draw node at (-1,-.5) [vertex] [] {};
	\draw node at (1,.5) [vertex] [] {};
	\draw node at (1,-.5) [vertex] [] {};
	\draw node at (-1,.5) [vertex] [] {};
	
	\draw[-,decorate] (-1,.5) -- (-1,-.5) node [midway, xshift=-.2cm ] {$e_1$};
	\draw[-] (0,0) -- (1,.5);
	\draw[-] (0,0) -- (1,-.5);
	\draw[-] (0,0) -- (-1,-.5);
	\draw[-] (0,0) -- (-1,.5);
	\draw[-,decorate] (1,.5) -- (1,-.5) node [midway, xshift =.2cm] {$e_2$};
	\end{tikzpicture}
	\caption{$e_1$ and $e_2$ are onbr adjacent}
	\label{onbr}
\end{figure}

Figure \ref{onbr} gives an example of two edges which are onbr adjacent. 

We can now introduce two more matchings based on cnbr and onbr. 

\begin{defn}
	A matching $M \subseteq E(G)$ is an onbr matching if the edges of $M$ are pairwise onbr independent. 
\end{defn}

\begin{defn}
A matching $M \subseteq E(G)$ is a cnbr matching if the edges of $M$ are pairwise cnbr independent.
\end{defn}

The parameters of interest for onbr and cnbr matchings are 
\[	\beta_{on} = \max \{ |M| : M \text{ is an onbr matching} \} \]
\[	\beta^-_{on} = \min \{ |M| : M \text{ is a maximal onbr matching} \} \]
\[	\beta_{cn} = \max \{ |M| : M \text{ is a cnbr matching} \} \]
\[	\beta^-_{cn} = \min \{ |M| : M \text{ is a maximal cnbr matching} \} \]

\section{Complexity Results}

Even a brief survey would be incomplete without mentioning some of the complexity results concerning matchings. Perhaps the most significant result is due to Edmonds \cite{Edmonds} celebrated blossom algorithm which shows that $\beta_1(G)$ can be determined in polynomial time for any $G.$ This result also has implications for some of the varieties of matchings discussed in Section 3. That is, since $\beta_1(G) = \beta_c(G) = \beta_{if}(G)$ if $G$ is connected, then $\beta_c(G)$ and $\beta_{if}(G)$ can also be determined in polynomial time for any $G.$ 

In contrast, Cameron \cite{Cameron} showed in 1989 that finding $\beta_*(G)$ is NP-complete for bipartite graphs. In 1992, Fricke and Laskar \cite{Fricke} showed finding $\beta_*(G)$ is linear for trees. Golumbic and Laskar \cite{Golumbic2} showed $\beta_*(G)$ can be found in polynomial time for interval graphs, chordal graphs, and circular arc graphs. Golumbic and Lewenstein \cite{Golumbic3} later generalized these results to wider classes of graphs, and improved on some of the previously known results. In particular, they showed that determining $\beta_*(G)$ is NP-complete for planar graphs, polynomial time for trapezoid graphs and cocomparability graphs, and linear time for interval graphs.

Golumbic, Lewenstein, and Hirst \cite{Golumbic1} showed that finding $\beta_{ur}(G)$ is NP-complete for bipartite graphs and chordal graphs. Furthermore, for proper interval graphs and threshold graphs, computational time is on the order of $\mathcal{O}(|V|).$ 
Moreover, this result can be extended to show that finding $\beta_{ac}(G)$ is also NP-complete. 

Of the matching varieties discussed in Section 3, this leaves only the complexity of $\beta_{dc}$ in question. 

We also note that for trees and for graphs with no even cycles, $\beta_1 = \beta_{ur}.$ That is, for graphs whose blocks are either edges or chordless odd cycles, we have $\beta_1 = \beta_{ur}.$

\section{Open Questions}

We end this survey with a few directions for further work concerning the $P$-matchings introduced here. 

Of course the first direction would be to find good upper bounds on $\beta_P$ and good lower bounds on $\beta^-_P$ for each of the $P$-matchings introduced. 

Also for each $P$-matching, characterization theorems should be developed to describe when a given matching $M$ is a maximum $P$-matching. 

Recalling the Gallai identities, it would be of interest to find a parameter, say ``$\alpha_P$" such that $\beta_P(G) + \alpha_P(G) = |V(G)|.$ 

It would also be desirable to have Nordhaus-Gaddum type results for these parameters, i.e., to complete the following inequalities
\[	? \leq \beta_P(G) + \beta_P(\bar{G}) \leq ?	\]
\[	? \leq \beta_P(G)\cdot \beta_P(\bar{G}) \leq ?	\]

In addition to many other results which could be of interest concerning these parameters, determining the complexity of finding $\beta_P$ is most desirable.


\begin{thebibliography}{1}
	\bibitem{Alavi} Alavi, Y., Behzad, M., Lesniak-Foster, L. M., \& Nordhaus, E. A. (1977). Total matchings and total coverings of graphs. \textit{Journal of Graph Theory,} 1(2), 135-140.
	
	\bibitem{Berge1} Berge, C., \& Minieka, E. (1973). Graphs and hypergraphs.
	
	\bibitem{Berge2} Berge, C. (1957). Two theorems in graph theory. \textit{Proceedings of the National Academy of Sciences,} 43(9), 842-844.
	
	\bibitem{Cameron} Cameron, K. (1989). Induced matchings. \textit{Discrete Applied Mathematics,} 24(1-3), 97-102.
	
	
	\bibitem{Edmonds} Edmonds, J. (1965). Paths, trees, and flowers. \textit{Canadian Journal of mathematics,} 17(3), 449-467.
	
	\bibitem{Fricke} Fricke, G., \& Laskar, R. (1992). Strong matchings on trees. \textit{Congressus Numerantium,} 239-239.
	
	\bibitem{Goddard} Goddard, W., Hedetniemi, S. M., Hedetniemi, S. T., \& Laskar, R. (2005). Generalized subgraph-restricted matchings in graphs. \textit{Discrete Mathematics,} 293(1-3), 129-138.
	
	\bibitem{Golumbic1} Golumbic, M. C., Hirst, T., \& Lewenstein, M. (2001). Uniquely restricted matchings. \textit{Algorithmica,} 31(2), 139-154.
	
	\bibitem{Golumbic2} Golumbic, M. C., \& Laskar, R. C. (1993). Irredundancy in circular arc graphs. \textit{Discrete Applied Mathematics,} 44(1-3), 79-89.
	
	\bibitem{Golumbic3} Golumbic, M. C., \& Lewenstein, M. (2000). New results on induced matchings. Discrete Applied Mathematics, 101(1-3), 157-165.
	
	\bibitem{Goodman} Goodman, S., Hedetniemi, S., \& Tarjan, R. E. (1976). B-matchings in trees. \textit{SIAM Journal on Computing,} 5(1), 104-108.
	
	\bibitem{Lovasz} Lov\'{a}sz, L., \& Plummer, M. D. (2009). \textit{Matching theory} (Vol. 367). American Mathematical Soc.
	
	\bibitem{Graham} Graham, R. L. (1970). On primitive graphs and optimal vertex assignments. \textit{Annals of the New York academy of sciences,} 175(1), 170-186.
	
	\bibitem{Mirsky} Mirsky, Leonid. \textit{Transversal theory.} Vol. 197. Academic Press, New York, 1971.
	
	\bibitem{Tutte} Tutte, W. T. (1947). The factorization of linear graphs. \textit{Journal of the London Mathematical Society,} 1(2), 107-111.
\end{thebibliography}
\end{document}